\documentclass[11pt,a4paper]{article}
\usepackage[latin1]{inputenc} 
\usepackage[T1]{fontenc}
\usepackage{amsmath, amssymb, dsfont,eufrak}
\usepackage{graphics, pstricks, epsfig, psfig}
\usepackage[frenchb]{babel}
\usepackage[all]{xy}
\usepackage{float}
\usepackage{multicol}
\usepackage{here}
\usepackage{array}

\def\color[#1]#2{}

\def\T0{\Theta_{0 0}}
\def\T1{\Theta_{0 1}}

\def\ep{\epsilon}

\def\epsilon{\varepsilon}

\def\theta{\vartheta}

\newcommand{\car}[2]{\left[ \begin{array}{c} #1 \\ #2 \end{array} \right]}

\newtheorem{theo}{Th{\'e}or{\`e}me}

\newtheorem{lemme}{Lemme}
         
\newenvironment{rem}{\textbf{Remarque} :}{}

\newtheorem{propo}{Proposition}

\title{
M{\'e}thode A.G.M. pour les courbes ordinaires de genre $3$
non hyperelliptiques sur $\mathbb{F}_{2^N}$}  
\author{Christophe Ritzenthaler}

\begin{document}
\maketitle

\begin{abstract}
We propose a A.G.M. algorithm for the determination of the
characteristic polynomial of an ordinary 
non hyperelliptic curve of genus $3$ over $\mathbb{F}_{2^N}$.
\end{abstract}
{\def\thefootnote{}\footnotetext{\hskip-1.8em Univ. Paris VII, laboratoire de
    th{\'e}orie des nombres. \\  E-mail : ritzenth@math.jussieu.fr}}%

Nous nous proposons dans la pr{\'e}sente note de donner un algorithme de
calcul du polyn{\^o}me caract{\'e}ristique de la jacobienne d'une courbe
de genre $3$ ordinaire et non-hyperelliptique sur
$k=\mathbb{F}_{2^N}$. Cet algorithme appartient {\`a} la cat{\'e}gorie des
m{\'e}thodes dites A.G.M. pour le calcul du nombre de points introduites
par J.-F. Mestre dans les cas hyperelliptiques (cf. \cite{mestre1} et
\cite{mestre2}). Nous illustrons son fonctionnement 
avec le calcul r{\'e}alis{\'e} sous MAGMA en $25$ heures pour une courbe sur
$\mathbb{F}_{2^{72}}$. Une impl{\'e}mentation plus efficace (utilisant les
m{\'e}thodes de
multiplication rapide) est en
cours. Les d{\'e}monstrations des r{\'e}sultats utilis{\'e}s seront
disponibles dans la th{\`e}se de l'auteur. Notons enfin que la m\'ethode
analogue dans le cas r\'eel permet sans doute le calcul d'une matrice de
Riemann de la jacobienne d'une courbe r\'eelle, imitant en cela
l'A.G.M. classique de Gauss dans le cas des courbes de genre $1$.

\section{M{\'e}thode A.G.M.} \label{para1}
Etant donn{\'e}e une courbe $\tilde{C}$
 de genre $g$ sur $k=\mathbb{F}_{2^N}$ ordinaire,
la mise en \oe uvre d'une m{\'e}thode A.G.M. pour la d{\'e}termination du
polyn{\^o}me caract{\'e}ristique proc{\`e}de globalement de la mani{\`e}re
suivante : 
\begin{enumerate}
\item on rel{\`e}ve la courbe sur une extension $K$  non
ramifi{\'e}e de $\mathbb{Q}_2$. Ce rel{\`e}vement ne doit pas {\^e}tre
quelconque si l'on veut pouvoir effectuer les calculs dans une
extension non ramifi{\'e}e. Nous donnerons dans le cas du genre $3$ au
paragraphe \ref{bonmodele}, un
bon mod{\`e}le sur une extension de degr{\'e} au plus $4 N$ de $\mathbb{Q}_2$.
\item Une fois en caract{\'e}ristique $0$, par analogie avec le cas
  complexe, 
on calcule alors en fonction
  des coefficients de ce mod{\`e}le
   $2^g$ rapports de 
<<th{\^e}ta constantes>> de la forme suivante
$$\left(\frac{\theta \car{f}{e}(0)}{\theta \car{0}{0}(0)}\right)^2_{f,e \in
  (\mathbb{Z}/2\mathbb{Z})^g}$$
dans le cas du genre $3$ les $8$
   rapports qui nous int{\'e}r{\`e}ssent
   parmi les $2^{g-1}(2^g+1)=36$ non nulles sont exactement ceux
   pour lesquelles le rapport pr{\'e}c{\'e}dent est congru {\`a} $1$ modulo
   $8$. Ces derniers orientent en effet le noyau de l'isog{\'e}nie pour
   que celle-ci coincide avec celui du rel{\`e}vement du Verschiebung qui est un
   sous-espace totalement isotrope pour le couplage de Weil. 
Quitte {\`a} faire alors un changement de base
   symplectique on peut supposer que ceux-ci sont tels que
   $f=0$. 
 On note  ces rapports $(A_e^{(0)})_{e \in (\mathbb{Z}/2\mathbb{Z})^g}.$\\   
Il faut exprimer ces rapports en fonction des coefficients de la courbe.
Dans le cas hyperelliptique, les
   relations sont donn{\'e}es par la
   formule de Thomae \cite[livre II]{mumford}. L'objet du paragraphe
   \ref{constante} est de donner un analogue de ces formules pour le
   cas du genre $3$ non hyperelliptique. Alors que dans
   le cas hyperelliptique  sont utilis{\'e}s des invariants binaires
   reli{\'e}s aux points de Weierstrass, on utilise ici des invariants
   ternaires bas{\'e}s sur des produits de d{\'e}terminants de bitangentes.
\item On duplique ces constantes par une g{\'e}n{\'e}ralisation de la
  moyenne arithm{\'e}tico-g{\'e}om{\'e}trique pour le genre $1$ : 
$$A_e^{(i+1)}=\frac{1}{2^g} \sum_{f \in (\mathbb{Z}/2
  \mathbb{Z})^g} A_e^{(i)} \sqrt{\frac{A_{e+f}^{(i)}}{A_e^{(i)}}}$$ 
la racine carr{\'e}e de $x \in 1+8 \mathcal{O}$ {\'e}tant choisie comme 
celle congrue
{\`a} $1$ modulo $4$.
\item On peut alors montrer que $
      \left(A_e^{(N(n+1))}/A_e^{(N)}\right)$  
converge vers $\alpha=\pm \pi_1\ldots\pi_g$ o{\`u} les $\pi_i$ sont les
$g$ racines du Frobenius inversibles modulo $2$ (\cite{mestre2}).
 Ce th{\'e}or{\`e}me peut
{\^e}tre d{\'e}montr{\'e} par la th{\'e}orie de la multiplication complexe, l'argument
de convergence provenant d'un th{\'e}or{\`e}me de R. Carls
(cf. \cite{carls}).
\item On calcule alors le polyn{\^o}me minimal de
  $\beta=\alpha+2^{gN}/\alpha$. Dans le cas $g=1,2$ on peut montrer que le
  polyn{\^o}me minimal est essentiellement donn{\'e} par ce nombre
  unique qu'il suffit de conna{\^\i}tre avec une pr{\'e}cision $N/2$
  respectivement $N$ ce qui implique respectivement $3N/2$ et $2 N$
  it{\'e}rations. Ce n'est plus le cas en genre $3$ : on utilise alors un
  analogue $2$-adique de la m{\'e}thode LLL pour la d{\'e}termination du
  polyn{\^o}me minimal. L'utilisation brutale 
de cette m{\'e}thode  semble
  toutefois n{\'e}cessiter la d{\'e}termination de $\alpha$ avec une
  pr{\'e}cision {\'e}gale {\`a} $11 N$.  
\item Lorsque $g \leq 3$ on montre alors que la connaissance de ce
  polyn{\^o}me minimal permet de d{\'e}terminer le polyn{\^o}me
  caract{\'e}ristique au signe pr{\`e}s. Ce n'est plus toujours vrai en
  genre sup{\'e}rieur (cf. \cite{mestre2}).
\item Il reste alors {\`a} lever l'ambiguit{\'e} sur le signe. Cela peut
  se faire  en effectuant un calcul dans la jacobienne
  (ce qui permet {\'e}galement de v{\'e}rifier l'algorithme)
  ou si l'on conna{\^\i}t des points rationnels d'ordre $4$.
\end{enumerate}

Les deux paragraphes suivants donnent dans le cas $g=3$ non
hyperelliptique les principaux {\'e}l{\'e}ments qui permettent de mettre
en \oe uvre la m{\'e}thode d{\'e}crite, en particulier la d{\'e}termination
des th{\^e}ta constantes et le choix d'un bon mod{\`e}le. Comme il a
{\'e}t{\'e} dit, les justifications th{\'e}oriques seront donn{\'e}es
ult{\'e}rieurement.

\section{D{\'e}termination des th{\^e}ta constantes} \label{constante}
 Soit $K$ est un corps alg{\'e}briquement clos de caract{\'e}ristique $0$ et
$C/K$ une courbe de genre $3$ non hyperelliptique que l'on consid{\`e}re
plong{\'e}e dans $\mathbb{P}^2$ comme une quartique lisse.
Nous utilisons ici des travaux de Riemann et Weber (cf. \cite{riemann}
et \cite{weber}). Ces
derniers sont bas{\'e}s sur le calcul des $28$ bitangentes {\`a} $C$ (i.e
les droites $l$ telles que $(l \cdot C)=2 P+2 Q$) 
que nous explicitons comme suit. \\

On note $\Sigma=\{L \in \textrm{Pic}^2(C), L^2=\mathcal{K}\}$
l'ensemble des fibr{\'e}s th{\^e}ta caract{\'e}ristiques (o{\`u} $\mathcal{K}$
d{\'e}signe le fibr{\'e} canonique) qui se partagent en $\Sigma_0$ et
$\Sigma_1$ sous-ensembles pour lesquels $\textrm{h}^0(L)=0$ (resp. $1$)
et dits  fibr{\'e}s pairs (resp. impairs). On montre alors que
l'ensemble des bitangentes est canoniquement en bijection avec
$\Sigma_1$.\\
D'autre part, le choix d'une base symplectique de $\textrm{Jac}(C)[2]$
permet d'associer {\`a} chaque point de $2$-torsion une matrice $2
\times 3$ not{\'e}e $\car{\ep}{\ep'} $ {\`a} coefficients dans $\mathbb{F}_2$
(appel{\'e}e caract{\'e}ristique) et de r{\'e}aliser une bijection explicite
avec $\Sigma$. On montre alors qu'{\`a} chaque bitangente est associ{\'e}e
une caract{\'e}ristique impaire (i.e. pour laquelle
$\ep \cdot \ep'=1$).\\
Afin de mieux comprendre la combinatoire de ces
caract{\'e}ristiques, on introduit un sous-ensemble de l'ensemble des
caract{\'e}ristiques appel{\'e} ensemble
principal not{\'e} $([i])_{i=1 \ldots 7}$ caract{\'e}ris{\'e} par les propri{\'e}t{\'e}s suivantes :
\begin{itemize}
\item Toute caract{\'e}ristique impaire s'{\'e}crit soit $[i]$ soit
  $[i]+[j]$, $i \neq j$.
\item Toute caract{\'e}ristique paire s'{\'e}crit soit $[0]$ soit
  $[i]+[j]+[k]$ avec $i,j,k$ distincts.\\
\end{itemize}

On a par exemple
$$\begin{array}{cccc}
[1] =\left[ \begin{array}{ccc} 1 & 0 & 1 \\ 0 & 0&1 \end{array}\right]
&
[2] =\left[ \begin{array}{ccc} 1 & 1 & 0 \\ 0 & 1&1 \end{array}\right]
&
[3] =\left[ \begin{array}{ccc} 1 & 1 & 1 \\ 0 & 1&0 \end{array}\right]
&
[4] =\left[ \begin{array}{ccc} 0 & 1 & 1 \\ 1 & 0&1 \end{array}\right]
\vspace{.2cm}
\\{}
[5] =\left[ \begin{array}{ccc} 1 & 0 & 0 \\ 1 & 0&0 \end{array}\right]
&
[6] =\left[ \begin{array}{ccc} 0 & 0 & 1 \\ 1 & 1&1 \end{array}\right]
&
[7] =\left[ \begin{array}{ccc} 0 & 1 & 0 \\ 1 & 1&0 \end{array}\right]
&
\end{array}$$
On note alors $\beta_i$ la bitangente associ{\'e}e {\`a} $[i]$ et
$\beta_{ij}$ celle associ{\'e}e {\`a} $[i]+[j]$.\\

Quitte {\`a} r{\'e}aliser une transformation lin{\'e}aire on peut supposer
$$\begin{cases}
\beta_1 : x_1=0 & \beta_5=a_1 x_1 +a_2 x_2 +a_3 x_3 \\
\beta_2 : x_2=0 & \beta_6=a_1' x_1+a_2' x_2 +a_3' x_3 \\
\beta_3 : x_3=0 & \beta_7=a_1'' x_1+a_2'' x_2 +a_3'' x_3 \\
\beta_4 : x_1+x_2+x_3=0 &
\end{cases}$$
on a alors le r{\'e}sultat suivant d{\^u} {\`a} Riemann \cite{riemann}
\begin{theo} \label{theo1}
Il existe un plongement de $C$ dans $\mathbb{P}^2$ de la forme
$$\sqrt{x_1 u_1}+\sqrt{x_2 u_2}+\sqrt{x_3 u_3}=0$$
o{\`u} $u_1,u_2,u_3$ sont lin{\'e}aires en $x_1,x_2,x_3$ donn{\'e}s par
$$\begin{cases}
u_1 + u_2 + u_3 +x_1 +x_2+x_3 =0 \\
\frac{u_1}{a_1} + \frac{u_2}{a_2} +\frac{u_3}{a_3}+k a_1 x_1 +k a_2
x_2+k a_3
x_3=0\\
\frac{u_1}{a_1'} + \frac{u_2}{a_2'} +\frac{u_3}{a_3'}+k' a_1' x_1 +k'
a_2' x_2+k' a_3'
x_3=0\\
\frac{u_1}{a_1''} + \frac{u_2}{a_2''} +\frac{u_3}{a_3''}+k'' a_1'' x_1
+k'' a_2''
x_2+k'' a_3'' 
x_3=0\\
\end{cases}$$
avec  $k,k',k''$ solutions des  deux syst{\`e}mes
lin{\'e}aires suivants :
$$\left( \begin{array}{ccc}
\frac{1}{a_1} & \frac{1}{a_1'} & \frac{1}{a_1''} \\
\frac{1}{a_2} & \frac{1}{a_2'} & \frac{1}{a_2''} \\
\frac{1}{a_3} & \frac{1}{a_3'} & \frac{1}{a_3''} \\
\end{array} \right) \left( \begin{array}{c} \lambda\vphantom{\frac{1}{a_3''}} \\ \lambda'\vphantom{\frac{1}{a_3''}} \\ \lambda''\vphantom{\frac{1}{a_3''}}
\end{array} \right) =\left( \begin{array}{c} -1 \vphantom{\frac{1}{a_3''}}\\ -1 \vphantom{\frac{1}{a_3''}} \\ -1 \vphantom{\frac{1}{a_3''}}
\end{array} \right),$$
$$
\left( \begin{array}{ccc}
\lambda a_1 & \lambda' a_1' & \lambda'' {a_1''} \\
\lambda {a_2} & \lambda' {a_2'} & \lambda'' {a_2''} \\
\lambda {a_3} & \lambda' {a_3'} & \lambda'' {a_3''} \\
\end{array} \right) \left( \begin{array}{c} k \\ k' \\ k''
\end{array} \right) =\left( \begin{array}{c} -1 \\ -1 \\ -1
\end{array} \right).
$$

Sous cette forme une expression des $28$ bitangentes est alors donn{\'e}e par
$$
\begin{array}{ccc}
\beta_1 : x_1=0 & \beta_2 : x_2=0 & \beta_3 : x_3=0 \\
\beta_{23} : u_1=0 & \beta_{13}  : u_2=0 & \beta_{12} : u_3=0 \\
\end{array}$$
$$\begin{array}{cc}
\beta_4 : x_1+x_2+x_3=0  & \beta_5 : a_1 x_1+a_2 x_2+a_3 x_3=0 \\ 
\beta_6 : a_1' x_1+a_2' x_2 +a_3' x_3=0 & \beta_7 : a_1'' x_1 +a_2'' x_2 +a_3''
x_3=0 \\
\end{array}$$
$$\begin{array}{cc}
\beta_{14}: u_1+x_2+x_3=0 & \beta_{15}: 
\frac{u_1}{a_1}+k a_2 x_2+k a_3 x_3=0 \\
\beta_{16} :\frac{u_1}{a_1'} + k' a_2' x_2+ k' a_3' x_3=0 & \beta_{17} : \frac{u_1}{a_1''} + k'' a_2'' x_2+ k'' a_3'' x_3=0
\end{array}$$
$$\begin{array}{cc}
\beta_{24} :x_1+u_2+x_3=0 & \beta_{25} : k a_1 x_1+\frac{u_2}{a_2}+k a_3 x_3=0 \\
\beta_{26} : k' a_1' x_1+\frac{u_2}{a_2'}+ k' a_3' x_3=0 & \beta_{27} :
k'' a_1'' x_1+ \frac{u_2}{a_2''} + k'' a_3'' x_3=0 
\end{array}$$
$$\begin{array}{cc}
\beta_{34} : x_1+x_2+u_3=0 & \beta_{35} : k a_1 x_1+k a_2 x_2+ \frac{u_3}{a_3}=0 \\
\beta_{36} : k' a_1' x_1+k' a_2' x_2+ \frac{u_3}{a_3'}=0 & \beta_{37} :
k'' a_1'' x_1+ k'' a_2'' x_2+\frac{u_3}{a_3''}=0
\end{array}$$
$$\begin{array}{c}
\beta_{67} : \frac{u_1}{1-k a_2 a_3} +\frac{u_2}{1- k a_3 a_1} +
\frac{u_3}{1-k a_1 a_2}=0 \\
\beta_{57} : \frac{u_1}{1-k' a_2' a_3'} +\frac{u_2}{1- k' a_3' a_1'} +
\frac{u_3}{1-k' a_1' a_2'}=0 \\
\beta_{56} : \frac{u_1}{1-k'' a_2'' a_3''} +\frac{u_2}{1- k'' a_3'' a_1''} +
\frac{u_3}{1-k'' a_1'' a_2''}=0 \\
\end{array}
$$
$$\begin{array}{c}
\beta_{45} : \frac{u_1}{a_1(1-k a_2 a_3)} +\frac{u_2}{a_2 (1- k a_3 a_1)} +
\frac{u_3}{a_3  (1-k a_1 a_2)}=0 \\
\beta_{46} : \frac{u_1}{a_1'(1-k' a_2' a_3')} +\frac{u_2}{a_2' (1- k' a_3'
a_1')}  +
\frac{u_3}{a_3'  (1-k' a_1' a_2')}=0 \\
\beta_{47} : \frac{u_1}{a_1''(1-k a_2'' a_3'')} +\frac{u_2}{a_2'' (1- k'' a_3''  a_1'')} +
\frac{u_3}{a_3''  (1-k'' a_1'' a_2'')}=0 \\
\end{array}$$
\end{theo}

Si on se donne $C$ sous une forme
$\sqrt{x_1 u_1}+\sqrt{x_2 u_2}+\sqrt{x_3 u_3}=0$
il nous reste {\`a} montrer comment d{\'e}terminer les 
$\beta_i$. On peut r{\'e}aliser cela au moyen du calcul suivant

\begin{enumerate}
\item On calcule $D_1(\lambda)$, d{\'e}terminant de la hessienne de la
  famille de coniques suivante
$$
Q_1(\lambda)=\lambda ^2 (x_2 u_3)+ \lambda (x_1 u_1-x_2 u_2-x_3 u_3) + (x_3
u_2).$$
\item  On calcule $R_1(x_1,x_2,x_3)=\textrm{Res}(D_1,Q_1,\lambda)$.
\item On calcule de m{\^e}me $R_2(x_1,x_2,x_3)$ avec la famille
  $$Q_2(\lambda)=\lambda ^2 (x_1 u_3)+ \lambda (x_2 u_2-x_1 u_1-x_3
  u_3) + (x_3 u_1).$$
\item Enfin $R=\textrm{pgcd}(R_1,R_2)/(x_3 u_3)= \prod_{i=4}^7 \beta_i$. 
\end{enumerate}

La d{\'e}termination des bitangentes est essentielle au calcul des
th{\^e}ta constantes. On a en effet le r{\'e}sultat suivant essentiellement
contenu dans Weber (cf. \cite{weber})
\begin{theo} \label{th2}
Soit $[\chi]=[i]+[j]+[k]$. Alors
$$\left( \frac{\vartheta[\chi](0)}{\vartheta(0)}
\right)^4=\frac{[\beta_{i},\beta_{j},\beta_{ij}] [\beta_{ik},\beta_{jk},\beta_{ij}] [\beta_{j},\beta_{jk},\beta_{k}] [\beta_{i},\beta_{ik},\beta_{k}]}{[\beta_{j},\beta_{jk},\beta_{ij}]
  [\beta_{i},\beta_{ik},\beta_{ij}]
[\beta_{i},\beta_{j},\beta_{k}] [\beta_{ik},\beta_{jk},\beta_{k}]}$$
o{\`u} l'on note $[\beta_{l_1},\beta_{l_2},\beta_{l_3}]=\det(\beta_{l_1},\beta_{l_2},\beta_{l_3})$. 
\end{theo}

\section{Bon mod{\`e}le} \label{bonmodele}
On reprend ici les notations du paragraphe \ref{para1}.
Une courbe de genre $3$ non-hyperelliptique sur
$k=\mathbb{F}_{2^N}$ a au plus $7$ bitangentes et
le mod{\`e}le de Riemann $\sqrt{x_1 u_1}+\sqrt{x_2 u_2}+\sqrt{x_3
  u_3}=0$  n'est plus valable en caract{\'e}ristique $2$. On commence
donc par chercher une description simple des courbes qui nous int{\'e}ressent
\begin{lemme}
Toute courbe de genre $3$ ordinaire, non hyperelliptique sur $k$
 est isomorphe   {\`a} une courbe d'{\'e}quation
\begin{eqnarray} \label{modele}
(a x^2+ b y^2+ c z^2+ d x y +e x z + f y z)^2-x y
z (x+y+z)=0
\end{eqnarray}
avec la condition suivante 
$$a b c (a+b+d) (a+c+e) (b+c+f) (a+b+c+d+e+f+1) \ne 0.$$ 
Inversement toutes les courbes qui v{\'e}rifient ces conditions sont des
courbes de genre $3$ ordinaires et non hyperelliptiques.\\
De plus l'isomorphisme est d{\'e}fini sur une extension au plus cubique
de $k$.
\end{lemme}

On consid{\`e}re donc $\tilde{C}$ donn{\'e}e par les coefficients $a,b,c,d,e,f$
ci-dessus. Nous voulons remonter ce mod{\`e}le sur une extension de
$\mathbb{Q}_2$ non ramifi{\'e}e de telle sorte que ses bitangentes
soient faciles {\`a} d{\'e}terminer, c'est-{\`a}-dire qu'elles soient
d{\'e}finies  par les formules du
th{\'e}or{\`e}me \ref{theo1}. Nous allons pour cela l{\'e}g{\'e}rement d{\'e}former
le mod{\`e}le de Riemann.\\ 

Tout comme la forme $y^2=x^3+a x+b$ doit {\^e}tre
transform{\'e}e en caract{\'e}ristique $2$ en $y^2+xy=x^3+ax+b$, nous
allons transformer notre mod{\`e}le de Riemann
 en rajoutant des termes \textit{{\`a} la}
Artin-Schreier. Remarquons pour cela que le mod{\`e}le de Riemann
$\sqrt{x_1 u_1}+\sqrt{x_2 u_2}+\sqrt{x_3 u_3}=0$ s'{\'e}crit {\'e}galement 
$$\begin{cases}
Y_1^2=x_1 u_1\\
Y_2^2=x_2 u_2\\
Y_3^2=x_3 u_3\\
Y_1+Y_2+Y_3=0
\end{cases}$$
Consid{\'e}rerons le mod{\`e}le suivant sur $k$
$$\begin{cases}
Y_1^2+l_1 Y_1=l_1 v_1 \\
Y_2^2+l_2 Y_2=l_2 v_2 \\
Y_3^2+l_3 Y_3=l_3 v_3 \\
l_1+l_2+l_3=0 \\
Y_1+Y_2+Y_3=l
\end{cases}$$
o{\`u} $l_1,l_1,l_2,v_1,v_2,v_3$ et $l$ sont lin{\'e}aires. On a par
exemple
\begin{propo} \label{propo1}
$\tilde{C}$ 
donn{\'e}e par le mod{\`e}le (\ref{modele}) est
 isomorphe sur $k$ {\`a} la courbe donn{\'e}e par le mod{\`e}le
ci-dessus avec 
$$\begin{cases}
l_1=x,l_2=y,l_3=x+y\\
l=z \\
v_1= b c y+(c+f) z \\
v_2= a c x + d c y + (c+f) z \\
v_3= a c x + (d+b) c y + (1+e+c+f) z
\end{cases}$$
\end{propo}
On rel{\`e}ve alors les expressions obtenues sur $K$ {\`a} l'identique sauf
 pour $l_3=-2 z-x-y$ et on {\'e}crit
$$\begin{cases}
(2 Y_1+l_1)^2=l_1(4 v_1+l_1) \\
(2 Y_2+l_2)^2=l_2(4 v_2+l_2) \\
(2 Y_3+l_3)^2=l_3(4 v_3+l_3) \\
(2 Y_1+u_1)+(2 Y_2+u_2)+(2 Y_3+u_3)=0
\end{cases}$$
Si on effectue le changement de coordonn{\'e}es
$x=x_1,y=x_2,z=-(x_1+x_2+x_3)/2$ on a le mod{\`e}le de Riemann 
\begin{equation} \label{mod2}
C  : \sqrt{x_1\vphantom{(l_1)} \smash{\underbrace{(4 v_1+l_1)}_{= u_1}}}
+\sqrt{x_2\vphantom{(l_1)} \smash{
  \underbrace{(4 v_2+l_2)}_{= u_2}}}+\sqrt{x_3 \vphantom{(l_1)} \smash{\underbrace{(4
  v_3+l_3)}_{=u_3}}}=0
\end{equation}

\section{Exemple}
On consid{\`e}re la courbe suivante :
$$\tilde{C}:
(x^2+\omega^2 y^2+\omega^4 z^2+\omega x y+(\omega^2+1) x
z+(\omega^3+1) y z)^2-xyz(x+y+z)=0$$ 
d{\'e}finie sur $k=F_{2^{N}}$, $N=72$ o{\`u} $\omega$ engendre le groupe
multiplicatif $k^*$.\\
Les calculs sont r{\'e}alis{\'e}s gr{\^a}ce au logiciel MAGMA version 2.9
sur un Pentium III {\`a} 1.13 Ghz avec 2 GigaOctets de m{\'e}moire centrale.\\
On calcule alors les $v_1,v_2,v_3$ du mod{\`e}le de la proposition \ref{propo1}
\begin{eqnarray*}
v_1 &=& \omega^6 y + (\omega^4 + \omega^3 + 1) z \\
v_2 &=&    \omega^4 x + \omega^5 y + (\omega^4 + \omega^2 + 1) z\\
v_3 &=&    \omega^4 x + (\omega^6 + \omega^5) y +  (\omega^4 + \omega^3 +
\omega^2 + 1) z
\end{eqnarray*}
On note  $K$ l'extension de degr{\'e} $d$ non ramifi{\'e}e de
$\mathbb{Q}_2$, $\mathcal{O}$ son anneau d'entiers et $w \in
\mathcal{O}$ 
qui se r{\'e}duit sur $\omega$.
Le mod{\`e}le (\ref{mod2}) de $C$ est d{\'e}fini par
\begin{eqnarray*}
u_1 &=&(-2 w^4 - 2 w^3 - 1) x_1 + (4 w^6 - 2 w^4 - 2 w^3 - 2) x_2 + (-2 w^4 - 
    2 w^3 - 2) x_3 \\
u_2 &=&
(2 w^4 - 2 w^2 - 2) x_1 + (4 w^5 - 2 w^4 - 2 w^2 - 1)x_2 + (-2 w^4 - 2
w^2 - 2) x_3 \\
 u_3&=&
(2 w^4 - 2 w^3 - 2 w^2 - 2) x_1 + (4 w^6 + 4 w^5 - 2 w^4 - 2 w^3 - 2 w^2 
    - 2) x_2 + \\ & & (-2 w^4 - 2 w^3 - 2 w^2 - 1) x_3
\end{eqnarray*}
On effectue  le calcul des $4$ bitangentes qui nous manquent gr{\^a}ce {\`a}
l'algorithme du paragraphe \ref{constante}
(temps : $1/2$ heure, elles sont d{\'e}finies sur $K$ dans notre cas)
 puis de l'ensemble des bitangentes et de tous les
rapports non nuls
 $$\left(\frac{\theta \car{f}{e}(0)}{\theta \car{0}{0}(0)}\right)^4_{f,e \in
  (\mathbb{Z}/2\mathbb{Z})^g}$$ gr{\^a}ce aux th{\'e}or{\`e}mes
\ref{theo1} et \ref{th2}. On s{\'e}lectionne  les $8$ rapports congrus {\`a} $1$ modulo
$16$ et on calcule leur racine carr{\'e}e dans $K$ (temps : $15$ mn)
$$\begin{cases}
 1 + O(2^4) \\
1 + (w^{16} + w^{14} + w^{12} + w^{11} + w^{9} + w^{5} + w^{4}) \cdot 2^3 +
O(2^4) \\
1 + (w^{16} + w^{14} + w^{12} + w^{11} + w^{10} + w^8 + w^4 + w^3)
\cdot 2^3 + O(2^4) \\
1 + (w^{10} + w^9 + w^8 + w^5 + w^3) \cdot 2^3 + O(2^4) \\
1 + (w^6 + w^4 + w^3) \cdot 2^3 + O(2^4) \\
1 + (w^{16} + w^{14} + w^{12} + w^{11} + w^9 + w^6 + w^5 + w^3) \cdot
2^3 + O(2^4) \\
 1 + (w^{16} + w^{14} + w^{12} + w^{11} + w^{10} + w^8 + w^6) \cdot 2^3 + O(2^4)
 \\
1 + (w^{10} + w^9 + w^8 + w^6 + w^5 + w^4) \cdot 2^3 + O(2^4)
\end{cases}$$

 On
effectue alors l'it{\'e}ration de la formule de duplication $12 N$
fois (temps : 24 heures) qui nous donne au signe pr{\`e}s 
la valeur du produit des racines
inversibles du Frobenius modulo $2$ 
$$\beta=1 + 2^8 + 2^9 + 2^{11} + 2^{13} + 2^{15} +2^{16} + 2^{17} +
2^{18} +
 2^{21}
+ 2^{23} + 2^{24} + \ldots + 2^{787} + 2^{790} + 2^{791} + O(2^{793})$$
On trouve alors gr{\^a}ce {\`a} LLL (temps : $1$ seconde) le polyn{\^o}me minimal
de $\beta$ :
\begin{eqnarray*}
P_{\textrm{sym}}(x) &=&
 x^4- 
52767044410803560460262696266497 x^3- \\
& &
 78121277277710794719527572033891108646286909 \cdot 2^{72} x^2+ \\
& & 610161746623391968394415142270976679056928051874538813 \cdot 2^{2
 \cdot 72} x +\\
& &46918330565326150855288775851644884890720289023905899509851903489
 \cdot 2^{3 \cdot 72} 
\end{eqnarray*}
D'o{\`u} le polyn{\^o}me caract{\'e}ristique au signe pr{\`e}s
\begin{eqnarray*}
P(x) &=&x^6- 9925657555 x^5+ 1108548370771462406931
x^4 \\
& & -146512229527151304651245280013057 x^3+
 1108548370771462406931 \cdot 2^{72} x^2 \\
& & - 9925657555 \cdot 2^{2 \cdot 72} x+2^{3
   \cdot 72}.
\end{eqnarray*}

Pour d{\'e}terminer le signe, on {\'e}tudie alors les points d'ordre $4$
de la jacobienne. Plus pr{\'e}cis{\'e}ment dans notre cas, $2^6$ divise $P(1)$
mais $16$ ne divise pas $P(-1)$. Il suffit donc que la jacobienne de
$C$ poss{\`e}de un point d'ordre $4$ rationnel pour conclure que $P(X)$
est le <<bon>> polyn{\^o}me. Consid{\'e}rons par exemple le point d'ordre
$2$ $(P_1+Q_1-(P_2+Q_2))=\frac{1}{2}((x \cdot C)-(y \cdot C))$. On
cherche un point de la jacobienne sous la forme $D=P+Q+R-3 P_2$ tel que
$2 D=P_1+Q_1-(P_2+Q_2)$ ou encore $2(P+Q+R)-(P_1+Q_1+5 P_2-Q_2)=0$.
En utilisant par exemple MAGMA, on montre que l'espace $L(P_1+Q_1+5
P_2-Q_2)$ qui est de dimension $4$, contient une fonction d\'efinie
sur $\mathbb{F}_{2^{72}}$ ayant trois z\'eros d'ordre $2$, d'o\`u
l'existence d'un point d'ordre $4$ rationnel sur
$\mathbb{F}_{2^{72}}$ (temps : $10$ secondes).\\

\begin{rem}
A posteriori on constate que la pr{\'e}cision n{\'e}cessaire est seulement
de $10 \cdot 72$, le calcul est alors r{\'e}alis{\'e} en $21$ heures. Il
serait donc tr{\`e}s utile d'avoir une meilleure borne pour LLL.
\end{rem}

\end{document}